\documentclass[a4paper,11pt,reqno]{amsart}
\usepackage{amsmath}
\usepackage{amstext}
\usepackage{amsbsy}
\usepackage{amsopn}
\usepackage{upref}
\usepackage{amsthm}
\usepackage{amsfonts}
\usepackage{amssymb}
\usepackage{mathrsfs}
\usepackage{times}
\allowdisplaybreaks
 \newtheorem{theorem}{Theorem}

 \newtheorem{lemma}[theorem]{Lemma}
\theoremstyle{definition}
 \newtheorem{definition}{Definition}
\theoremstyle{remark}
 
\newcommand{\p}{\partial}
\newcommand{\sgn}{\operatorname{sgn}}
\newcommand{\supp}{\operatorname{supp}}
\newcommand{\re}{\operatorname{Re}}
\newcommand{\im}{\operatorname{Im}}
\begin{document}
\title[dispersive equation]{The initial value problem for a third order dispersive equation on the two-dimensional torus}
\author[H.~Chihara]{Hiroyuki CHIHARA}
\address{Mathematical Institute,  
         Tohoku University, 
         Sendai 980-8578, Japan}
\email{chihara@math.tohoku.ac.jp}
\begin{abstract}
We present the necessary and sufficient condition for 
the $L^2$-well-posedness of the initial problem 
for a third order linear dispersive equation 
on the two dimensional torus. 
Birkhoff's method of asymptotic solutions is 
used to prove the necessity. 
Some properties of a system for 
quadratic algebraic equations associated to 
the principal symbol 
play crucial role in proving the sufficiency.  
\end{abstract}
\maketitle
\section{Introduction}
\label{section:1}
This paper is concerned with the initial value problem of the form 
\begin{alignat}{2}
Lu&=f(t,x)
&
\quad\text{in}\quad
&
\mathbb{R}\times\mathbb{T}^2,
\label{equation:pde}
\\
u(0,x)&=u_0(x)
&
\quad\text{in}\quad
&
\mathbb{T}^2,
\label{equation:data}
\end{alignat}
where 
$u(t,x)$ is a real-valued unknown function of 
$(t,x)=(t,x_1,x_2)\in\mathbb{R}\times\mathbb{T}^2$, 
$\mathbb{T}^2=\mathbb{R}^2/2\pi\mathbb{Z}^2$, 
$u_0(x)$ and $f(t,x)$ are given real-valued functions, 
$$
L
=
\p_t+p(\p)
+
\sum_{\lvert\alpha\rvert=2}
\frac{2!}{\alpha!}
a_{\sigma(\alpha)}(x)\p^\alpha
+
\vec{b}(x)\cdot\nabla
+
c(x),
$$
$\p_t=\p/\p{t}$, 
$\p_j=\p/\p{x_j}$, 
$\p=\nabla=(\p_1,\p_2)$, 
$p(\xi)=\xi_1\xi_2(\xi_1+\xi_2)$, 
$\alpha=(\alpha_1,\alpha_2)$ is a multi-index, 
$\lvert\alpha\rvert=\alpha_1+\alpha_2$, 
$\alpha!=\alpha_1!\alpha_2!$, 
$\sigma(\alpha)=(\alpha_1-\alpha_2)/2$, 
$\vec{b}(x)=(b_1(x),b_2(x))$, 
and 
$a_{\sigma(\alpha)}$, $b_j(x)$, $c(x)$ 
are real-valued smooth functions on $\mathbb{T}^2$. 
Such operators arise 
in the study of gravity wave of deep water. 
See \cite{bks}, \cite{dysthe}, \cite{hogan} 
and the references therein. 
\par
The purpose of this paper is to present 
the necessary and sufficient condition of 
the existence of a unique solution to 
\eqref{equation:pde}-\eqref{equation:data}. 
To state our results, 
we introduce notation and the definition of $L^2$-well-posedness.  
$C^\infty(\mathbb{T}^2)$ 
is the set of all smooth functions on
$\mathbb{T}^2$. 
$L^2(\mathbb{T}^2)$ 
is the set of all square-integrable functions on
$\mathbb{T}^2$. 
For $g{\in}L^2(\mathbb{T}^2)$, set 
$$
\lVert{g}\rVert
=
\left(
\int_{\mathbb{T}^2}\lvert{g(x)}\rvert^2dx
\right)^{1/2}. 
$$
$C(\mathbb{R};L^2(\mathbb{T}^2))$ 
is the set of all $L^2(\mathbb{T}^2)$-valued 
continuous functions on $\mathbb{R}$, 
and similarly,  
$L^1_{\text{loc}}(\mathbb{R};L^2(\mathbb{T}^2))$ 
is the set of all $L^2(\mathbb{T}^2)$-valued 
locally integrable functions on $\mathbb{R}$. 
Here we state the definition of $L^2$-well-posedness. 
\begin{definition}
\label{definition:l2wp}
The initial value problem 
{\rm \eqref{equation:pde}-\eqref{equation:data}} 
is said to be $L^2$-well-posed 
if for any $u_0{\in}L^2(\mathbb{T}^2)$ 
and $f{\in}L^1_{\text{loc}}(\mathbb{R};L^2(\mathbb{T}^2))$, 
{\rm \eqref{equation:pde}-\eqref{equation:data}} 
possesses a unique solution 
$u{\in}C(\mathbb{R};L^2(\mathbb{T}^2))$. 
\end{definition}
In view of Banach's closed graph theorem, 
if \eqref{equation:pde}-\eqref{equation:data} is $L^2$-well-posed, 
then for any $T>0$, there exists $C_T>0$ such that 
all the solutions satisfy the energy inequality 
\begin{equation}
\lVert{u(t)}\rVert
\leqslant
C_T
\left(
\lVert{u_0}\rVert
+
\left\lvert
\int_0^t\lVert{f(s)}\rVert{ds}
\right\rvert
\right),
\quad
(\lvert{t}\rvert\leqslant{T}). 
\label{equation:energy}
\end{equation}
\par
Our results are the following.
\begin{theorem}
\label{theorem:main}
The following conditions are mutually equivalent. 
\\
{\rm I}. \quad
{\rm \eqref{equation:pde}-\eqref{equation:data}} is $L^2$-well-posed. 
\\
{\rm II}. \quad
For any $x\in\mathbb{T}^2$ and for any $\xi\in\Lambda$, 
\begin{equation}
 \int_0^{2\pi}
 \sum_{\lvert\alpha\rvert=2}
 \frac{2!}{\alpha!}
 a_{\sigma(\alpha)}(x+tp^\prime(\xi))\xi^\alpha
 dt
 =0, 
\label{equation:vanish}
\end{equation}
where 
$
p^\prime(\xi)=\nabla_\xi{p(\xi)}
=(\xi_2(2\xi_1+\xi_2),\xi_1(2\xi_2+\xi_1))$ 
and 
$
\Lambda
=
\{
\xi\in\mathbb{R}^2 
\vert 
p^\prime(\xi)\in\mathbb{Z}^2
\}
$. 
\\
{\rm III}. \quad
$a_0(x)=a_1(x)+a_{-1}(x)$,
and there exists $\phi(x){\in}C^\infty(\mathbb{T}^2)$ 
such that 
$\nabla\phi(x)=(a_{-1}(x),a_1(x))$.
\end{theorem}
Here we explain the background of our problem. 
There are many papers dealing with the well-posedness of 
the initial value problem for dispersive equations. 
Generally speaking, 
it is difficult to characterize the well-posedness. 
In fact, results on the characterization are very limited. 
In \cite{mizohata} Mizohata studied 
the initial value problem for 
Schr\"odinger-type operator of the form 
$$
S
=
\p_t-iq(\p)+\vec{b}(x)\cdot\nabla+c(x)
$$
on $\mathbb{R}^n$, 
where $i=\sqrt{-1}$ and $q(\xi)=\xi_1^2+\dotsb+\xi_n^2$. 
He gave the necessary condition for 
the $L^2$-well-posedness
\begin{equation}
\sup_{(T,x,\xi)\in\mathbb{R}\times\mathbb{R}^n\times\mathbb{R}^n}
\left\lvert
\int_0^T
\im\vec{b}(x+tq^\prime(\xi))\cdot\xi
dt
\right\rvert
<+\infty. 
\label{equation:mizohata}
\end{equation}
\eqref{equation:mizohata} 
gives the upper bound of the admissible bad first order term 
$(\im\vec{b}(x))\cdot\nabla$. 
In other words, 
\eqref{equation:mizohata} 
is necessary so that $(\im\vec{b}(x))\cdot\nabla$ 
can be controled by the local smoothing effect of 
$e^{itq(\p)}$. 
In \cite{ichinose1} 
Ichinose generalized \eqref{equation:mizohata} 
for the Schr\"odinger-type equation on a 
complete Riemannian manifold. 
When $n=1$, 
\eqref{equation:mizohata} 
is also the sufficient condition for $L^2$-well-posedness. 
Moreover, in \cite{tarama1} and \cite{tarama2} 
Tarama characterized the $L^2$-well-posedness 
for third order equations on $\mathbb{R}$. 
\par
For equations on compact manifolds, 
the local smoothing effect breaks down 
everywhere in the manifold. 
Then, the restriction on the bad lower order terms 
becomes stronger, 
and the characterization of $L^2$-well-posedness 
seems to be relatively easier. 
In fact, 
the $L^2$-well-posedness for $S$ on $\mathbb{T}^n$ is characterized. 
See \cite{chihara}, \cite{ichinose2} and \cite{tsukamoto}. 
\par
$L$ is the simplest example of 
higher order dispersive operators on 
higher dimensional spaces. 
$p(\xi)$ satisfies 
$p^\prime(\xi)\ne0$ and 
$\det{p^{\prime\prime}(\xi)}\ne0$ 
for $\xi\ne0$. 
The symbol of the Laplacian $q(\xi)$ also satisfies 
the same conditions. 
It is interesting that our method of the proof 
does not work if we replace $p(\xi)$ 
by $r(\xi)=\xi_1^3+\xi_2^3$. 
This seems to be due to the degeneracy of 
$\det{r^{\prime\prime}}(\xi)=36\xi_1\xi_2$ for $\xi\ne0$. 
\par
The organization of this paper is as follows. 
In Sections~\ref{section:2} and \ref{section:3} 
we prove I $\Longrightarrow$ II and 
II $\Longrightarrow$ III respectively. 
We omit the proof of III $\Longrightarrow$ I. 
In deed, under the condition 
$a_0=a_{1}+a_{-1}$, $L$ becomes 
$$
L
=
\p_t
+
p(\p)
+
a_{-1}(x)\frac{\p p}{\p\xi_1}(\p)
+
a_{1}(x)\frac{\p p}{\p\xi_2}(\p)
+
\vec{b}(x)\cdot\nabla
+
c(x).
$$
In view of $\nabla\phi=(a_{-1},a_1)$, 
we have 
$$
e^\phi{L}e^{-\phi}
=
\p_t+A, 
\quad
A
=
p(\p)
+
\tilde{b}_1(x)\p_1
+
\tilde{b}_2(x)\p_2
+
\tilde{c}(x),
$$
where 
$\tilde{b}_1(x)$, 
$\tilde{b}_2(x)$ 
and  
$\tilde{c}(x)$ 
are real-valued smooth functions on $\mathbb{T}^2$. 
It is easy to see that the initial value problem for 
$e^\phi{L}e^{-\phi}$ is $L^2$-well-posed since 
$$
A+A^\star
=
-\frac{\p\tilde{b}_1}{\p{x_1}}(x)
-\frac{\p\tilde{b}_2}{\p{x_2}}(x). 
$$
\section{The Proof of I $\Longrightarrow$ II}
\label{section:2}
We begin with the reduction of \eqref{equation:vanish}. 
Set 
$$
a(x,\xi)
=
\sum_{\lvert\alpha\rvert=2}
\frac{2!}{\alpha!}
a_{\sigma(\alpha)}(x)\xi^\alpha
$$
for short. If 
\begin{equation}
\sup_{(T,x,\xi)\in\mathbb{R}\times\mathbb{T}^2\times\mathbb{R}^2}
\left\lvert
\int_0^T
a(x+tp^\prime(\xi),\xi)
dt
\right\rvert
<+\infty, 
\label{equation:real}
\end{equation}
then \eqref{equation:vanish} holds 
since $a(x+tp^\prime(\xi),\xi)$ is a $2\pi$-periodic function in
$t\in\mathbb{R}$ for any $(x,\xi)\in\mathbb{T}^2\times\Lambda$. 
\eqref{equation:real} is equivalent to 
an apparently weaker condition 
\begin{equation}
\sup_{(T,x,\xi)\in\mathbb{R}\times\mathbb{T}^2\times\mathbb{Q}^2}
\left\lvert
\int_0^T
a(x+tp^\prime(\xi),\xi)
dt
\right\rvert
<+\infty. 
\label{equation:rational}
\end{equation}
For $\xi\in\mathbb{Q}^2$, there exist 
$\alpha\in\mathbb{Z}^2$ and $l\in\mathbb{N}$ such that 
$\xi=\alpha/l$. 
Changing the variable by $t=l^2s$, we have 
$$
\int_0^T
a(x+tp^\prime(\xi),\xi)
dt
=
\int_0^{T/l^2}
a(x+sp^\prime(\alpha),\alpha)
dt.
$$
Then, \eqref{equation:rational} is equivalent to 
\begin{equation}
\sup_{(T,x,\alpha)\in\mathbb{R}\times\mathbb{T}^2\times\mathbb{Z}^2}
\left\lvert
\int_0^T
a(x+tp^\prime(\alpha),\alpha)
dt
\right\rvert
<+\infty. 
\label{equation:integer}
\end{equation}
To prove I $\Longrightarrow$ II. 
we shall show the contraposition of 
I $\Longrightarrow$ \eqref{equation:integer}. 
\par
Suppose that \eqref{equation:integer} fails to hold, 
that is, 
for any $n\in\mathbb{N}$, 
there exist 
$T_n\in\mathbb{R}$, 
$x_n\in\mathbb{T}^2$ 
and 
$\alpha\in\mathbb{Z}^2$ 
such that 
$$
\left\lvert
\int_0^{T_n}
a(x_n+tp^\prime(\alpha),\alpha)
dt
\right\rvert
\geqslant
2n. 
$$
We shall show that the energy inequality 
\eqref{equation:energy} fails to hold. 
Firstly, we consider the case that $T_n>0$ and 
$$
\int_0^{T_n}
a(x_n+tp^\prime(\alpha),\alpha)
dt
\geqslant
2n. 
$$
Since $[0,T_n]\times\mathbb{T}^2$ is compact and 
$$
(s,x)\in[0,T_n]\times\mathbb{T}^2
\longmapsto
\int_0^s
a(x+tp^\prime(\alpha),\alpha)
dt
$$
is continuous, 
there exist $(T,y)\in(0,T_n]\times\mathbb{T}^2$ such that 
$$
\max_{(s,x)\in[0,T_n]\times\mathbb{T}^2}
\int_0^s
a(x+tp^\prime(\alpha),\alpha)
dt
=
\int_0^T
a(y+tp^\prime(\alpha),\alpha)
dt. 
$$
Pick up $\psi{\in}C^\infty(\mathbb{T}^2)$ such that 
$\lVert\psi\rVert=1$ and 
$$
\int_0^T
a(x+tp^\prime(\alpha),\alpha)
dt
\geqslant
n
\quad\text{in}\quad
\supp[\psi].
$$
\par
Let $u$ be a complex-valued solution to 
\eqref{equation:pde} 
with a complex-valued given function $f(t,x)$. 
Then, $L\re{u}=\re{f}$ and $L\im{u}=\im{f}$ 
since all the coefficients in $L$ are real-valued. 
We construct a sequence of complex-valued 
asymptotic solutions to $Lu=0$. 
For $l\in\mathbb{N}$, set 
$$
u_l(t,x)
=
e^{itp(l\alpha)+il\alpha\cdot{x}+\phi_l(t,x)}
\psi(x+(t-T/l^2)p^\prime(l\alpha)),
$$
$$
\phi_l(t,x)
=
\int_0^{l^2t}
a(x+sp^\prime(\alpha),\alpha)
ds.  
$$
Then, $u_l{\in}C^\infty(\mathbb{R}\times\mathbb{T}^2)$, 
\begin{align}
  \lVert{u_l}(0)\rVert
& =
  \lVert\psi\rVert
  =1,
\label{equation:initial}
\\
  \lVert{u_l}(T/l^2)\rVert
& =
  \lVert\exp(\phi_l(T/l^2,\cdot))\psi(\cdot)\rVert
  \geqslant
  e^n.
\label{equation:evolution}
\end{align}
\par
Next we compute $Lu_l$. 
Set $b(x,\xi)=\vec{b}(x)\cdot\xi$ and 
$v_l(t,x)=e^{-itp(l\alpha)-il\alpha\cdot{x}}u_l(t,x)$ 
for short. 
We deduce 
\begin{align}
  e^{-itp(l\alpha)-il\alpha\cdot{x}}Lu_l
& =
  (\p_t+ip(l\alpha))v_l
  +
  p(\p+il\alpha)v_l
\nonumber
\\
& +
  a(x,\p+il\alpha)v_l
  +
  b(x,\p+il\alpha)v_l
  +
  c(x)v_l,
\label{equation:lu} 
\\
  (\p_t+ip(l\alpha))v_l
& =
  ip(l\alpha)v_l+a(x+tp^\prime(l\alpha),l\alpha))v_l
\nonumber
\\
& +
  e^{\phi_l}
  p^\prime(l\alpha)\cdot\nabla\psi(x+(t-T/l^2)p^\prime(l\alpha))
\label{equation:part1}
\\
  p(\p+il\alpha)v_l
& =
  (p(\p)+il\alpha\cdot{p^\prime(\p)}-p^\prime(l\alpha)\cdot\nabla-ip(l\alpha))v_l
\nonumber
\\
& =
  (p(\p)+il\alpha\cdot{p^\prime(\p)})v_l-ip(l\alpha)v_l
\nonumber
\\
& -
  \left(
  \int_0^{l^2t}
  p^\prime(l\alpha)\cdot\nabla(a(x+sp^\prime(\alpha),\alpha))
  ds
  \right)
  v_l
\nonumber
\\
& -
  e^{\phi_l}
  p^\prime(l\alpha)\cdot\nabla\psi(x+(t-T/l^2)p^\prime(l\alpha))
\nonumber
\\
& =
  (p(\p)+il\alpha\cdot{p^\prime(\p)})v_l-ip(l\alpha)v_l
\nonumber
\\
& -
  \left(
  \int_0^{l^2t}
  \frac{d}{ds}(a(x+sp^\prime(\alpha),l\alpha))
  ds
  \right)
  v_l
\nonumber
\\
& -
  e^{\phi_l}
  p^\prime(l\alpha)\cdot\nabla\psi(x+(t-T/l^2)p^\prime(l\alpha))
\nonumber
\\
& =
  (p(\p)+il\alpha\cdot{p^\prime(\p)})v_l-ip(l\alpha)v_l
\nonumber
\\
& -
  a(x+tp^\prime(l\alpha),l\alpha))v_l
  +
  a(x,l\alpha)v_l
\nonumber
\\
& -
  e^{\phi_l}
  p^\prime(l\alpha)\cdot\nabla\psi(x+(t-T/l^2)p^\prime(l\alpha)),
\label{equation:part2}
\\
  a(x,\p+il\alpha)v_l
& =
  a(x,\p)v_l-a(x,l\alpha)v_l
\nonumber
\\
& +
  2il(\alpha_1a_1(x)+\alpha_2a_0(x))\p_1v_l
\nonumber
\\
& +
  2il(\alpha_1a_0(x)+\alpha_2a_{-1}(x))\p_2v_l,
\label{equation:part3}
\\
  b(x,\p+il\alpha)v_l
& =
  b(x,\p)v_l-ib(x,l\alpha)v_l.
\label{equation:part4} 
\end{align}
Substituting 
\eqref{equation:part1}, 
\eqref{equation:part2}, 
\eqref{equation:part3} 
and 
\eqref{equation:part4} 
into
\eqref{equation:lu}, 
we obtain 
$$
\lvert{Lu_l(t,x)}\rvert
\leqslant
Cl\lvert\alpha\rvert
\sum_{\lvert\beta\rvert\leqslant3}
\lvert\p^\beta{v_l(t,x)}\rvert.
$$
Then, we deduce 
\begin{equation}
\lvert{Lu_l(t,x)}\rvert
\leqslant
C_0l\lvert\alpha\rvert
(1+T\lvert\alpha\rvert^2)^3
\exp\left(
    \int_0^T
    a(y+sp^\prime(\alpha),\alpha)
    ds
    \right) 
\label{equation:part5}
\end{equation}
for $t\in[0,T/l^2]$. 
Integrating the $L^2(\mathbb{T}^2)$-norm of 
\eqref{equation:part5} over $[0,T/l^2]$, 
we have 
$$
\int_0^{T/l^2}
\lVert{Lu_l(t)}\rVert
dt
\leqslant
\frac{A_\alpha}{l},
$$
$$
A_\alpha
=
2\pi{C_0T}\lvert\alpha\rvert
(1+T\lvert\alpha\rvert^2)^3
\exp\left(
    \int_0^T
    a(y+sp^\prime(\alpha),\alpha)
    ds
    \right). 
$$
If we take $l$ satisfying $A_\alpha\leqslant{l}$, then 
\begin{equation}
\int_0^{T/l^2}
\lVert{Lu_l(t)}\rVert
dt
\leqslant
1.
\label{equation:rhs}
\end{equation}
Combining 
\eqref{equation:initial}, 
\eqref{equation:evolution} 
and 
\eqref{equation:rhs}, 
we obtain 
$$
\lVert{u_l(T/l^2)}\rVert
\geqslant
e^n
>
2
\geqslant
\lVert{u_l(0)}\rVert
+
\int_0^{T/l^2}
\lVert{Lu_l(t)}\rVert
dt,
$$
which breaks the energy inequality \eqref{equation:energy}.
\par
When $T_n>0$ and 
$$
\int_0^{T_n}
a(x_n+tp^\prime(\alpha),\alpha)
dt
\leqslant
-2n,  
$$
we employ a sequence of asymptotic solutions of the form 
$$
u_l(t,x)
=
e^{-itp(l\alpha)-il\alpha\cdot{x}-\phi_l(t,x)}
\psi(x+(t-T/l^2)p^\prime(l\alpha)).
$$
When $T_n<0$, the proof above works also in $[T,0]$ 
for some $T\in[T_n,0]$. 
The proof of I $\Longrightarrow$ II finished. 
\section{The Proof of II $\Longrightarrow$ III}
\label{section:3}
To prove II $\Longrightarrow$ III, 
we need to know the properties of $\Lambda$. 
\begin{lemma}
\label{theorem:alg}
For any $\alpha\in\mathbb{Z}^2$, 
there exists $\xi(\alpha)\in\Lambda$ such that 
$p^\prime(\pm\xi(\alpha))=\alpha$. 
Moreover, 
$\xi(0)=0$, 
and for $\alpha\ne0$, $\xi(\alpha)\ne0$ and 
\begin{equation}
\xi_1(\alpha)\xi_2(\alpha)
+
\xi_1(-\alpha)\xi_2(-\alpha)
\ne0. 
\label{equation:trick}
\end{equation}
\end{lemma}
\begin{proof}
For the sake of intelligibility, 
we express two-vectors by the entries as 
$(\xi,\eta)\in\Lambda$ 
and 
$(\alpha,\beta)\in\mathbb{Z}^2$. 
We solve a system for quadratic algebraic equations 
$$
\eta(2\xi+\eta)=\alpha, 
\quad
\xi(2\eta+\xi)=\beta.
$$
{\bf Case} $(\alpha,\beta)=(0,0)$.\quad
Suppose $\eta(2\xi+\eta)=0$ and $\xi(2\eta+\xi)=0$. 
Then $\eta=0$ or $2\xi+\eta=0$, 
and $\xi=0$ or $2\eta+\xi=0$. 
In any case, 
$(\xi,\eta)=0$ is a unique solution. 
\vspace{6pt}\\
{\bf Case} $\alpha=0$, $\beta\ne0$.\quad
Suppose $\beta\ne0$, 
$\eta(2\xi+\eta)=0$ and $\xi(2\eta+\xi)=\beta$. 
Then, $\eta=0$ or $\eta=-2\xi$, 
and $\xi(2\eta+\xi)=\beta$. 
If $\eta=0$, then $\xi^2=\beta$, 
which implies 
$\beta>0$ and $(\xi,\eta)=(\pm\sqrt{\beta},0)$. 
If $\eta=-2\xi$, then $-3\xi^2=\beta$, 
which implies 
$\beta<0$ and 
$(\xi,\eta)=\pm(\sqrt{-\beta/3},-2\sqrt{-\beta/3})$. 
Then, we have 
$$
(\xi(0,\beta),\eta(0,\beta))
=
\begin{cases}
\pm(\sqrt{\beta},0) & (\text{if}\quad\beta>0)
\\
\pm\left(\sqrt{-\dfrac{\beta}{3}},-2\sqrt{-\dfrac{\beta}{3}}\right) & (\text{if}\quad\beta<0) 
\end{cases}
$$ 
$$
\xi(0,\beta)\eta(0,\beta)
+
\xi(0,-\beta)\eta(0,-\beta)
=
-\frac{2}{3}\lvert\beta\rvert\ne0.
$$
\vspace{6pt}\\
{\bf Case} $\alpha\ne0$, $\beta=0$.\quad
In the same way as the case $\alpha=0$ and $\beta\ne0$, 
we have 
$$
(\xi(\alpha,0),\eta(\alpha,0))
=
\begin{cases}
\pm(0,\sqrt{\alpha}) & (\text{if}\quad\alpha>0)
\\
\pm\left(-2\sqrt{-\dfrac{\alpha}{3}},\sqrt{-\dfrac{\alpha}{3}}\right) & (\text{if}\quad\alpha<0) 
\end{cases}
$$ 
$$
\xi(\alpha,0)\eta(\alpha,0)
+
\xi(\alpha,0)\eta(0,-\alpha)
=
-\frac{2}{3}\lvert\alpha\rvert\ne0.
$$
\vspace{6pt}\\
{\bf Case} $\alpha\beta\ne0$.\quad
Suppose  $\alpha\beta\ne0$, 
$\eta(2\xi+\eta)=\alpha$ and $\xi(2\eta+\xi)=\beta$. 
$\xi\eta\ne0$ 
since 
$\alpha\beta=\xi\eta(2\eta+\xi)(2\xi+\eta)\ne0$. 
Substituting 
$\eta=-\xi/2+\beta/2\xi$ 
into 
$\eta(2\xi+\eta)=\beta$, 
we have 
$3\xi^4+2(2\alpha-\beta)\xi^2-\beta^2=0$. 
Then, 
$$
\xi^2
=
\frac{\beta-2\alpha\pm\sqrt{(\beta-2\alpha)^2+3\beta^2}}{3}.
$$
Since 
$\lvert\beta-2\alpha\rvert<\sqrt{(\beta-2\alpha)^2+3\beta^2}$ 
and $\xi^2>0$, 
$$
\xi^2
=
\frac{\beta-2\alpha+\sqrt{(\beta-2\alpha)^2+3\beta^2}}{3}. 
$$
Then, we have 
$$
\xi(\alpha,\beta)
=
\pm
\sqrt{
\frac{\beta-2\alpha+\sqrt{(\beta-2\alpha)^2+3\beta^2}}{3}
}. 
$$
Using $2\xi\eta=\beta-\xi^2=\alpha-\eta^2$, we get 
\begin{align}
  2\xi\eta
& =
  \frac{2(\alpha+\beta)-\sqrt{(\beta-2\alpha)^2+3\beta^2}}{3}, 
\label{equation:xieta}
\\
  \eta^2
& =
  \frac{\alpha-2\beta+\sqrt{(\beta-2\alpha)^2+3\beta^2}}{3}. 
\label{equation:eta}
\end{align}
Here we remark that 
$\eta^2>0$ is satisfied in \eqref{equation:eta} 
since 
$$
(\sqrt{(\beta-2\alpha)^2+3\beta^2})^2-(\alpha-2\beta)^2
=
3\alpha^2>0.
$$
Using \eqref{equation:xieta} and \eqref{equation:eta}, 
we deduce 
\begin{align*}
  2\lvert\xi\eta\rvert
& =
  \chi(\alpha,\beta)
  \frac{2(\alpha+\beta)-\sqrt{(\beta-2\alpha)^2+3\beta^2}}{3},
\\
  \chi(\alpha,\beta)
& =
  \sgn\left(2(\alpha+\beta)-\sqrt{(\beta-2\alpha)^2+3\beta^2}\right).
\end{align*}
$\chi(\alpha,\beta)$ makes sense 
for $\alpha\beta\ne0$ since 
$$
4(\alpha+\beta)^2
-
\left(\sqrt{(\beta-2\alpha)^2+3\beta^2}\right)^2
=
12\alpha\beta\ne0.
$$
Thus, we have 
$$
\begin{bmatrix}
\xi(\alpha,\beta)
\\
\eta(\alpha,\beta)
\end{bmatrix}
=
\pm
\begin{bmatrix}
\sqrt{
\dfrac{\beta-2\alpha+\sqrt{(\beta-2\alpha)^2+3\beta^2}}{3}
} 
\\
\chi(\alpha,\beta)
\sqrt{
\dfrac{\alpha-2\beta+\sqrt{(\beta-2\alpha)^2+3\beta^2}}{3}
}
\end{bmatrix}.
$$
It follows from \eqref{equation:xieta} that 
$$
\xi(\alpha,\beta)\eta(\alpha,\beta)
+
\xi(-\alpha,-\beta)\eta(-\alpha,-\beta)
=
-\frac{1}{3}
\sqrt{(\beta-2\alpha)^2+3\beta^2}
\ne0.
$$
This completes the proof.
\end{proof}
Finally, we prove II $\Longrightarrow$ III. 
Express $a_{\sigma(\alpha)}(x)$ by the Fourier series of the form
$$
a_{\sigma(\alpha)}(x)
=
\sum_{\beta\in\mathbb{Z}^2}
a_{\sigma(\alpha),\beta}
\exp(i\beta\cdot{x}). 
$$
Substitute the Fourier series into 
\eqref{equation:vanish}. 
Then, for $(x,\xi)\in\mathbb{T}^2\times\Lambda$, 
\begin{align*}
  0
& =
  \sum_{\beta\in\mathbb{Z}^2}
  \sum_{\lvert\alpha\rvert=2}
  \frac{2!}{\alpha!}
  a_{\sigma(\alpha),\beta}\xi^\alpha
  e^{i\beta\cdot{x}}
  \int_0^{2\pi}
  e^{it\beta\cdot{p^\prime(\xi)}}
  dt
\\
& =
  2\pi
  \sum_{\substack{\beta\in\mathbb{Z}^2
                  \\
                  \beta\cdot{p^\prime(\xi)}=0}}
  (a_{1,\beta}\xi_1^2
   +
   2a_{0,\beta}\xi_1\xi_2
   +
   a_{-1,\beta}\xi_2^2)
  e^{i\beta\cdot{x}}
\\
& =
  2\pi
  \sum_{\substack{\beta\in\mathbb{Z}^2
                  \\
                  \beta\cdot{p^\prime(\xi)}=0}}
  \{
  (a_{-1,\beta},a_{1,\beta})\cdot{p^\prime(\xi)}
  +
  2(a_{0,\beta}-a_{1,\beta}-a_{-1,\beta})\xi_1\xi_2
  \}
  e^{i\beta\cdot{x}}.
\end{align*}
It follows that if
$\beta\cdot{p^\prime(\xi)}=0$ 
and 
$\xi\in\Lambda$, 
then 
$$ 
(a_{-1,\beta},a_{1,\beta})\cdot{p^\prime(\xi)}
+
2(a_{0,\beta}-a_{1,\beta}-a_{-1,\beta})\xi_1\xi_2
=0. 
$$
In view of Lemma~\ref{theorem:alg}, we have 
\begin{align}
  (a_{-1,\beta},a_{1,\beta})\cdot\alpha
  +
  2(a_{0,\beta}-a_{1,\beta}-a_{-1,\beta})
  \xi_1(\alpha)\xi_2(\alpha)
& = 0, 
\label{equation:part11}
\\
  -
  (a_{-1,\beta},a_{1,\beta})\cdot\alpha
  +
  2(a_{0,\beta}-a_{1,\beta}-a_{-1,\beta})
  \xi_1(-\alpha)\xi_2(-\alpha)
& = 0, 
\label{equation:part12}
\end{align}
for $\alpha\in\mathbb{Z}^2$ 
satisfying $\alpha\cdot\beta=0$. 
The sum of \eqref{equation:part11} 
and \eqref{equation:part12} is 
$$
2(a_{0,\beta}-a_{1,\beta}-a_{-1,\beta})
(\xi_1(\alpha)\xi_2(\alpha)
 +
 \xi_1(-\alpha)\xi_2(-\alpha))
=0.
$$
In view of \eqref{equation:trick}, 
we get $a_{0,\beta}=a_{1,\beta}+a_{-1,\beta}$ 
for all $\beta\in\mathbb{Z}^2$. 
Thus, $a_0(x)=a_1(x)+a_{-1}(x)$, and 
\eqref{equation:part11} becomes 
\begin{equation}
(a_{-1,\beta},a_{1,\beta})\cdot\alpha=0
\quad\text{if}\quad
\alpha\cdot\beta=0.
\label{equation:part13}
\end{equation}
$(a_{-1,0},a_{1,0})=0$ 
since $\alpha\cdot0=0$ 
for all $\alpha\in\mathbb{Z}^2$. 
For $\beta\ne0$, \eqref{equation:part13} implies that 
there exists $\phi_\beta\in\mathbb{C}$ such that 
$(a_{-1,\beta},a_{1,\beta})=i\phi_\beta\beta$. 
If we set 
$$
\phi(x)
=
\sum_{\beta\ne0}
\phi_\beta
e^{i\beta\cdot{x}},
$$
then 
$$
\nabla\phi(x)
=
\sum_{\beta\ne0}
i\phi_\beta\beta
e^{i\beta\cdot{x}}
=
\sum_{\beta\ne0}
(a_{-1,\beta},a_{1,\beta})
e^{i\beta\cdot{x}}
=
(a_{-1}(x),a_1(x)),
$$
which is desired. 
The proof of II $\Longrightarrow$ III finished. 
\bibliographystyle{amsplain}

\begin{thebibliography}{10}
\bibitem{bks}M.~Ben-Artzi, H.~Koch and J.-C.~Saut, 
{\it Dispersion estimates for third order equations in two dimensions}, 
Comm.\ Partial Differential Equations 
{\bf 28} (2003), 1943--1974. 
\bibitem{chihara}H.~Chihara, 
{\it The initial value problem for Schr\"odinger equations on the torus}. 
Int.\ Math.\ Res.\ Not.\ {\bf 2002}, 789--820.
\bibitem{dysthe}K.~B.~Dysthe, 
{\it Note on a modification to the nonlinear Schr\"odinger equation for application to deep water}, 
Proc.\ Roy.\ Soc.\ London Ser.\ A {\bf 369} (1979), 
105--114.
\bibitem{hogan}S.~J.~Hogan, 
{\it The fourth-order evolution equation for deep-water gravity-capillary waves}, 
Proc.\ Roy.\ Soc.\ London Ser.\ A {\bf 402} (1985), 
359--372.
\bibitem{ichinose1}W.~Ichinose, 
{\it On $L^2$ well posedness of the Cauchy problem 
for Schr\"odinger type equations 
on the Riemannian manifold and the Maslov theory}, 
Duke Math.\ J. {\bf 56} (1988), 549--588.
\bibitem{ichinose2}\bysame, 
{\it A note on the Cauchy problem for Schr\"odinger type equations on the Riemannian manifold}, 
Math.\ Japon.\ {\bf 35} (1990), 205--213.
\bibitem{mizohata}S.~Mizohata, 
``On the Cauchy problem'', 
Notes and Reports in Mathematics in Science and Engineering {\bf 3}, 
Academic Press, Inc., Orlando, FL; 
Science Press, Beijing, 1985.
\bibitem{tarama1}S.~Tarama, 
{\it On the wellposed Cauchy problem for some dispersive equations}, 
J. Math.\ Soc.\ Japan {\bf 47} (1995), 143--158.
\bibitem{tarama2}\bysame, 
{\it Remarks on $L^2$-wellposed Cauchy problem for some dispersive equations},  
J. Math.\ Kyoto Univ.\ {\bf 37} (1997), 757--765.
\bibitem{tsukamoto}C.~Tsukamoto, 
{\it Integrability of infinitesimal Zoll deformations}, 
Geometry of geodesics and related topics (Tokyo, 1982), 
97--104, Adv.\ Stud.\ Pure Math.\ {\bf 3}, 
North-Holland, Amsterdam, 1984.
\end{thebibliography}

\end{document}